\documentstyle[11pt]{article}

\newtheorem{th}{Theorem}[section]
\newtheorem{de}[th]{Definition}

\newtheorem{lem}[th]{Lemma}
\newtheorem{co}[th]{Corollary}
\newtheorem{re}[th]{Remark}

{

\begin{document}

\large
\title {\bf GENERALIZED PROJECTION OPERATORS IN BANACH
SPACES:\\ PROPERTIES  AND APPLICATIONS}
\author{\bf Ya. I. Alber\\
Department of Mathematics\\
Technioin-Israel Institute of Technology\\
Haifa 32000, Israel}
\date{}
\maketitle

\section{Introduction}
\setcounter{equation}{0}

Metric projection \hskip 0.2cm operators  \hskip 0.2cm in Hilbert \hskip 0.2cm
and Banach spaces
are widely used in
different areas of mathematics such as
functional and numerical analysis, theory of optimization
and approximation and for  problems of optimal control
and operations research, nonlinear and stochastic programming and game
theory.

Metric projection operators  can be defined in a similar way
in Hilbert and Banach spaces. At the same time, they differ
significantly in their properties \cite{h,li}.

A metric projection operator in Hilbert space is a monotone
(accretive) and nonexpansive operator. It provides an absolutely best
approximation for arbitrary elements from Hilbert space by the
elements of convex closed sets .
This leads to a variety of applications of this operator
for investigating theoretical questions in analysis and for
approximation
methods.  (For details see \cite{breg,pet,zd,aus,ct,al1,al0}).

Metric projection operators in Banach space do not have the properties
mentioned above and their applications are not straightforward.
(See \cite{r,aln3}).

On the other hand, questions of smoothness and differentiability
of metric projection operators in Banach spaces were actively
investigated \cite{bj,li,aln2,al2}. New results in this field are
immediately used
in various applications. For example, recently established
in \cite{aln2,al2} properties of uniform continuity of these operators
were used in \cite{al3,al4} to prove  stability  of  the  penalty  and
quasisolution methods.

Two of  the most important applications of the  metric projection
operators in Hilbert spaces are as follows:
\begin{itemize}
\item  solve a variational inequality by the iterative-projection  method,
\item  find common point of convex sets by the iterative-projection method.
\end{itemize}
------------------------\\
This research was supported in part by the Ministry of Science Grant
3481-1-91 and by the Ministry of Absorption Center for Absorption in Science.

In Banach space these problems can not be solved in the framework of metric
projection operators. Therefore,
in the present paper we introduce  new generalized projection operators
in Banach space as a natural  generalization of  metric projection
operators in Hilbert space.  To demonstrate our approach, we apply
these operators for solving two problems mentioned above in Banach
space.

In Section 2 and Section 3 we introduce notations and recall some results from
t
theory of variational inequalities and theory of approximation. Then
in Section 4 and Section 5 we describe the properties of metric projection
operators  $P_\Omega$ in
Hilbert and Banach spaces and also formulate equivalence theorems
between variational inequalities and direct projection equations  with these
operators.
In Section 6 we discuss the existence of strongly
unique
 best approximations  based on Clarkson's and parallelogram  inequalities.

In Section 7 we introduce generalized projection operator $\Pi_\Omega$ which
acts from Banach space $B$ on  the convex closed set $\Omega$ in the
same
space $B$. Then we state its properties and give the convergence theorem
for the
method of successive generalized projections  used to find a common
point of
convex sets.

In Section 8 new generalized projection operator $\pi_\Omega$
acting from  conjugate Banach space $B^{*}$    on convex closed set
 $\Omega$ in  the space $B$ and its properties are examined.
Then a theorem of equivalence of
the solutions of variational inequality and operator equation with
operator $\pi_\Omega$ is presented.  It constitutes the basis for
construction of iterative-projection methods for nonlinear problems in
Banach space
(including smooth and nonsmooth optimization problems)

Lastly, in Section 9 we establish a connection between variational inequalities
and Wiener-Hopf equations in Banach spaces by means of metric and
generalized projection operators.

Note that the main properties of metric and generalized projection operators
in Banach spaces have been obtained by using principally new technique
including Banach space geometry, parallelogram inequalities, nonstandard
Lyapunov functionals and estimates of moduli of monotonicity and continuity
for duality mappings.

\section{Variational Inequalities. Problems of Optimization}
\setcounter{equation}{0}
In this section we recall some of the results from the
theory of variational inequalities and formulate a problem on the
equivalence
between solutions of the variational inequalities and corresponding
operator
equations.  These equations are solved by the iterative-projection
methods. This
yields an aproximation of solutions of the initial variational
inequalities.

Let $B$ be a real (reflexive) uniformly convex and uniformly smooth
Banach
space, $B^{*}$ its conjugate (dual) space, $ ||\cdot ||, ||\cdot
||_{B^{*}},
||\cdot ||_H$ norms in the Banach spaces $ B, B^{*}$ and in Hilbert
space $ H$.
 As usually we denote $<\varphi,x>$ a dual  product in B.
 This determines pairing between $\varphi  \in B^{*}$ and $x \in  B$
\cite{br}.
Let $\Omega $ be a  nonempty convex closed set in B.

\begin {de}  \label{gpo}
The operator $P_{\Omega}: B \rightarrow {\Omega} \subset B$
is called metric projection operator if it yields
the correspondence between an arbitrary point  $x \in  B$ and nearest
point
 $\bar x \in \Omega $ according to minimization problem
\begin{equation} \label{f1}
P_{\Omega}x = \bar x; \;\;\; \bar  x: ||x - \bar x|| =
\inf_{\xi \in \Omega} ||x - \xi||.
\end{equation}
 \end {de}

Under our conditions operator $P_{\Omega}$ is defined at any point $x
\in  B$ and
it is single-valued, i.e. there exists for each $x \in  B$ a unique
$projection $
$ \bar x$  called best approximation \cite{h}.

Let $A$ be an arbitrary operator acting from Hilbert space  H
to H ,
$\alpha$ an arbitrary fixed positive number  and ($\varphi,x$)
an inner product in H.  Let also $f \in H$.
It is well  known that (see, for instance,  \cite{cea})
\begin {th} \label{th1}
 The point $x \in \Omega  \subset  H$ is a solution of variational
inequality
\begin {equation} \label{f2}
(Ax-f, \xi - x) \ge 0, {\qquad }        \forall \xi \in\Omega
\end{equation}
 if and only if x is a solution of the following operator equation in H
\begin{equation} \label {f3}
x=P_{\Omega}(x - \alpha (Ax-f)).
\end{equation}
\end{th}

This is an important statement, because it provides a basis for
constructing
approximate (iterative) methods  in Hilbert spaces. The simplest  method
of this type can be described as follows
\begin{equation} \label{f4}
x_{n + 1} = P_{\Omega}(x_{n}-\alpha_n (Ax_{n} - f)), \qquad n=0,1,2...
\end{equation}
Under suitable choice of relaxation parameters $\alpha_n$ , one can
prove
that iterative process (\ref {f4}) converges strongly to the solution of
the
variational inequality (\ref {f2}).
It can be done for operator $A$ which have different structures and
different types of smoothness  \cite{al1,al0,cea,aus,dunn}.
Moreover, one can establish both stability and nonasymptotic estimates
of convergence rate \cite{al1,al0}.

     We want to  emphasize that  the problem of solving operator
equation
$Ax=f$ and the problem of minimization of the functional $u(x)$ on
$\Omega$ are
 realized as variational inequalities (\ref {f2}) for $\Omega =  H$ and
$$Ax=\partial u(x),{\quad } x \in\Omega, {\quad } f = 0 $$
respectively. Here $\partial u(x)$  is gradient  or subgradient of the
functional
 $u(x)$.

Now we consider more general and more complicated case
of the variational inequality
\begin {equation} \label{f5}
<Ax - f, \xi - x> \ge 0, {\qquad }       \forall \xi \in\Omega
\end{equation}
in Banach space $B$ with operator $A$ acting from  $B$ to $B^{\ast}$
\cite{aln3}.
There is a natural problem to formulate and to prove  an analogue of
Theorem \ref{th1} in Banach space, and then to use it as a basis to
construct
iterative-projection methods similar to (\ref{f4}).

It is quite obvious that  the Banach space analogue of the equation
(\ref{f3}) has the
following form
\begin {equation} \label{f6}
x=\Gamma_{\Omega}(Fx-\alpha(Ax-f))
\end {equation}
with operator F  acting from  B to $B^{\ast}$. The equation (\ref{f6})
is unusual  because operator ${\Gamma}_\Omega$  "projects" elements from
the dual space $B^{\ast}$ on the set $\Omega \subset B$.
Here one can not use metric projection
operator $P_{\Omega}$ for this purpose because it acts from  B
to B.  It turned out that a natural generalization of metric
projection operator in  Hilbert space leads to a new operator which
we call generalized projection operator:
 $$\pi_{\Omega}  : B^{\ast} \rightarrow  \Omega \subset B .$$
This automatically yields the following form of the equation (\ref{f6})
 $$x= \pi_{\Omega } (Jx-\alpha(Ax-f))$$
where $J: B\rightarrow B^{\ast}$ is a normalized duality mapping in $B$
\cite{br}.
The operator $J$ is one of the most significant operators  in nonlinear
functional
analysis. In particular, it is used in the theory of optimization and in
the theory of
monotone and accretive operators in Banach spaces.  It is determined by
the
expression
 $$ <Jx,x> = ||Jx||_{B^{*}}||x|| = ||x||^2.$$
 Note also that a duality mapping exists in each Banach space.
 In what follows we recall from \cite{al3} some of the examples of this
mapping
in spaces
$l^p, L^p $ , $ W^p_m$ , $\infty >p>1 $ :
\begin{itemize}
\item  (i) $ l^p:  Jx = ||x||^{2-p}_{l^p} y\in l^q, {\quad } x =
\lbrace  x_1, x_2,...\rbrace  $,
 $ y =   \lbrace  x_1 ||x_1||^{p-2}, x_2 ||x_2||^{p-2},...\rbrace, \\
 p^{-1} + q^{-1} = 1 $
\item  (ii) $ L^p:  Jx = ||x||^{2-p}_{L^p} |x|^{p-2}x \in L^q $
 \item  (iii) $ W^p_m:  Jx = ||x||^{2-p}_{W^p_m} \sum_
 {|\alpha | \leq m} (-1)^ {|\alpha
|}D^\alpha
   (|D^\alpha x |^{p-2} D^\alpha x ) \in W^q_{-m} $
 \end{itemize}
 Note that in Hilbert space  $J$ is  an identity operator.

 Now we  define the iterative method similar (\ref{f4}) as follows
\begin {equation} \label{k6}
x_{n + 1} = \pi_{\Omega}(Jx_{n}-\alpha_n (Ax_{n} - f)), {\qquad
}n=0,1,2...
\end{equation}
We will give full mathematical foundation for this method including
three basic
aspects: convergence, stability and estimates of convergence rate, in
forthcoming
paper.

\section{ Problems of Approximation. Common Points of Convex Sets}
\setcounter{equation}{0}
Second important problem which is investigated in this paper  using
projection operators can be formulated as follows:  find  common point
of an
ordered collection of convex and closed  (i.e. Chebyshev) sets  $\lbrace
\Omega_{1},
\Omega_{2 },...,\Omega_{ m} \rbrace $ in uniformly  convex Banach  space
$B$.
Here we assume that  sets $\lbrace  \Omega_{i},  i=1,2,...,m \rbrace$
have nonempty intersection  $\Omega _{*}= \bigcap _{i=1}^{m}
\Omega_{i}$.
Let us define a composition
\begin {equation} \label{f9}
P =  P_{1} \circ  P_{2} \;\circ \cdotp \cdotp \cdotp \circ \;
 P_{m},  \qquad     P_{i} = P_{\Omega_{i}}
\end{equation}
 and introduce  method of successive projections according to a formula
\begin {equation} \label{f10}
x_{n+1} = P^{n+1}x_0,  {\qquad }     n=0,1,2,...,  {\qquad } x_o \in B.
\end{equation}

Convergence of the iterative process  (\ref {f9})  and  (\ref {f10}) as
well as of similar
processes to the  point  $ x_{*} \in \Omega _{*} $ was proved before
only in Hilbert
space .  (See \cite{vn,breg,bc,ct,d,gpr}).  In the formulae (\ref {f9})
and
 (\ref {f10})
which describe method of successive projections one can use metric
projection
operators $P_{\Omega}$ in Banach space. However, up to this date
no proof was suggested for the convergence of (\ref {f9}) and (\ref
{f10})
in Banach space. The reason is that in Hilbert space H the metric
projection
operator satisfies the following  significant inequality
\begin {equation} \label{f3.3}
||P_\Omega{ x} - \xi||_H \leq||x-\xi||_H , {\qquad }        \forall \xi
\in\Omega
\end{equation}
which can be obtained  from the property of nonexpansiveness of  this
operator
in H
 \begin {equation} \label{f3.4}
 {||P_\Omega{x} - P_\Omega{y}||_H} \leq {||x- y||_H}.
\end{equation}
It satisfies also a much stronger property (see Section 4)
\begin {equation} \label{f13}
{||P_\Omega{ x} - x||^2}_H
\leq {||x-\xi||^2}_H - {||P_\Omega{ x} - \xi||^2}_H ,  {\qquad }
\forall \xi \in\Omega .
\end{equation}
But in Banach spaces these properties do not hold in general
\cite{sm,r,gr}.

For example, in \cite{sm}  it is shown that in uniformly convex  Banach
space
with modulus of convexity   $\delta (\epsilon) $ of order  $\epsilon ^q
$,
$q\ge 2 $ \cite{dis},  the inequality
\begin {equation} \label{f14}
{||P_\Omega{ x} - x||}^q \leq {||x-\xi||}^q - \lambda  ||P_\Omega{ x} -
\xi||^q , {\qquad }
        \forall \xi \in\Omega .
 \end{equation}
holds. Coefficient $\lambda < 1$ in (\ref{f14}) and it depends on $q$.
Namely, in \cite{sm} it is defined in Banach spaces of the type
$ L^p and  W^p_m,\; 1 < p < \infty $ as follows (cf. (\ref {f161}) and (\ref
{f1
$$1 < p \leq 2,\;\; q=2, \;  \;  \lambda = (p-1)/8;$$
and
$$ 2 < p < \infty, \;\; q = p , \; \;  \lambda = 1/p2^p.$$

Inequality (\ref{f14}) yields
\begin {equation}  \label{f15}
{||P_\Omega{ x} - \xi||}^q \leq {\lambda}^ {-1}{||x-\xi||}^q -
 {\lambda}^ {-1}{||P_\Omega{ x} - x||}^q,
{\qquad }  \forall \xi \in\Omega.
\end{equation}
This does not guarantee the nonexpansiveness of   the metric projection
operator in
Banach space even for $\xi = y\in\Omega$  while in (\ref {f3.4}) $\xi =
y $ is
an arbitrary element of the space H. But without this property,
one can not to study the method (\ref {f9}), (\ref {f10}).

Now  we consider more general case.
In \cite{sm} a strongly unique best approximation  was defined as
follows.
\begin{de} \label{def1}
$\bar x$ is called a strongly unique best approximation
in $\Omega$ for the element  $x  \in  B$ if there exists a constant
$\lambda$
and a strictly increasing function $\phi (t) : R^{+}  \rightarrow R^{+}$
such that
 $\phi(0)=0$ and
\begin {equation}  \label {f151}
{\phi(||P_\Omega{ x} - x||)} \leq {\phi (||x-\xi||)} - \lambda \phi
(||P_\Omega{ x} - \xi||) ,
{\qquad }  \forall \xi \in\Omega.
\end{equation}
\end{de}

The projection $\bar x = P_\Omega  x$ in Hilbert space, and the
projection
$\bar x$ in Banach space under the conditions of  \cite{sm}, are
strongly unique best
 approximations in $\Omega $.

We  call the projection $\bar x $ {\it absolutely best approximation \/}
of $x\in B$
with respect  to  function $\phi (t) $ if  $\lambda  = 1$ in (\ref
{f151}).
In this case the inequality (\ref {f151}) can be represented in the
equivalent form
 $$\phi (||P_\Omega{ x} - \xi||) \leq {\phi (||x-\xi||)} - \phi
(||P_\Omega{ x} - x||),
 {\qquad }  \forall \xi \in\Omega. $$
It is clear from (\ref {f13}) that metric projection $\bar x$ in Hilbert space
is absolutely
best approximation with respect to  function  $\phi (t) = t^2 $ (or with
respect to
functional $\phi (\xi) =  {||x-\xi||^2}_H $ with  $x$  fixed). But it is not
true in Banach spaces.

Thus, metric projection operator can not be used in (\ref {f9}) and
(\ref {f10}).
Instead  we introduce new generalized projection operator
$$\Pi_{\Omega}:  B \rightarrow {\Omega} \in B $$
so that the inequalities (\ref{f3.3}) and (\ref{f13}) hold for  some
Lyapunov
functional $\phi (\cdot)$. In Section 7 we will provide the convergence
theorem for the process
(\ref {f9}) and (\ref {f10}) which now has the form
\begin {equation} \label{f16}
x_{n+1} = \Pi^{n+1}x_0,  {\qquad }     n=0,1,2,..., {\qquad } x_o \in B
\end{equation}
and
\begin {equation}  \label{f17}
\Pi =  \Pi_{1} \circ  \Pi_{2} \;\circ \cdotp \cdotp \cdotp \circ \;
 \Pi_{m},  \qquad      \Pi_{i} = \Pi_{\Omega_{i}}.
 \end{equation}

\section{ Metric Projection Operator  $P_{\Omega}$ in Hilbert Space}
\setcounter{equation}{0}
All results described in Section 2 and Section 3 for two basic problems were
obtained
only in Hilbert space. This is due to the fact that many remarkable
properties
of the metric projection operators can not be extended from Hilbert
space to
Banach space. This is why we introduce in Section 7 and Section 8 new
generalize
projection operators in Banach spaces which have all properties of
metric projection operators in Hilbert space .

Before that in Section 4 and Section 5 we compare complete lists of the
properti
of the
metric  projection operators in Hilbert  and Banach spaces.

We denote $ \bar x = P_{\Omega}x $. Let $\xi \in {\Omega}$ and $
{\Omega}
\subset H$.
The following properties are  valid in Hilbert space:
\cite {aus,breg,cea,gr,h,pet,al3,na,za}.

{\bf 4.a.}  $P_{\Omega}$ is fixed at each point $\xi $,
 $P_{\Omega}\xi = \xi.$

{\bf 4.b.}  $P_{\Omega}$ is monotone (accretive) in $ H$, i.e. for all
$x,y \in H $
$$(\bar x - \bar y, x-y) \ge 0. $$

{\bf 4.c.} $ (x - \bar x, \bar x -\xi) \ge 0,  {\qquad }        \forall
\xi \in\Omega  . $ \\

{\bf 4.d.} $ (x - \xi, \bar x -\xi) \ge 0,  {\qquad }        \forall \xi
\in\Omega . $

{\bf 4.e. } $ (x - \bar x,  x -\xi) \ge 0,  {\qquad }        \forall \xi
\in\Omega. $\\
In fact, even stronger inequality
$$ (x - \bar x,  x -\xi) \ge  ||x - \bar x||^2_H ,  \;\; \forall \xi
\in\Omega$$
holds  (see (\ref {f52}).

{\bf 4.f.}   $P_{\Omega}$ is nonexpansive in $ H$,  i.e.

 $${ ||\bar x - \bar y||_H} \leq {||x- y||_H} .$$

{\bf 4.g.}  $P_{\Omega}$ is $P$-strongly monotone in $ H$,  i.e.
 $$(\bar x - \bar y,  x-y) \ge  ||\bar x - \bar y||^2_H .$$

{\bf 4.h. } The operator $P_{\Omega}$ yields an absolutely best
approximation
of $x\in H$ with respect to the functional  $ {V_1 (x,\xi )=  ||x -
\xi||^2 _H}$

$${||\bar x  - \xi||}^2 _H \leq {||x-\xi||}^2 _H - {||x - \bar x||}^2_H
,
{\qquad }
\forall \xi \in\Omega .$$

{\bf 4.i.} Any  $P_{\Omega}$ satisfies the inequality
$$ ((I - P_{\Omega})x - (I - P_{\Omega})y, P_{\Omega}x - P_{\Omega}y) \ge 0,
\;\; \forall x,y \in H.$$

\section{ Metric Projection Operator  $P_{\Omega}$ in  Banach Space}
\setcounter{equation}{0}
Here we show that some of the properties of the metric projection
operators in Hilbert
space are not satisfied in Banach space. At the same time, we describe
in detail the
properties  of uniform continuity of the metric projection operators in
Banach space.

We denote $ \bar x = P_{\Omega}x $. Let $\xi \in {\Omega}$ and $
{\Omega}\subset B$.
The following properties hold in Banach space (the sign "-" from
5.a - 5.i denotes an absence of corresponding property):

{\bf  5.a.} The operator $P_{\Omega}$ is fixed at each point $\xi $, i.e.
$P_{\Omega}\xi = \xi.$

{\bf  5.b.} -

{\bf  5.c.}  $ < J(x - \bar x), \bar x -\xi > \ge 0,  {\qquad }
\forall \xi \in\Omega  $ (see \cite{li}).

{\bf  5.d.} -

{\bf  5.e. } $ < J(x - \bar x),  x -\xi > \ge 0,  {\qquad } \forall \xi
\in\Omega $
(see \cite{al3}).\\
In what follows we show that a stronger statement  is true \cite{al3}.
\begin {th} \label{f51}
$\bar x \in \Omega $ is a projection of the point $x \in B$ on $\Omega$
if and only
if the inequality
\begin {equation} \label{f52}
  < J(x - \bar x),  x -\xi > \ge ||x - \bar x||^2 ,  {\qquad }  \forall
\xi \in\Omega
\end {equation}
is satisfied.
\end {th}

In fact, from (\ref {f52}) it follows immediately that
$$||x - \bar x|| \leq  {||x - \bar x||}^{ -1} { < J(x - \bar x),  x -\xi
> } \leq  ||x - \xi||,
\qquad   \forall \xi \in \Omega $$
 i.e.  $\bar x = P_{\Omega}x $.  Inversely, if   $ \bar x = P_{\Omega}x
$,  then by
virtue of  5.c  we have
\begin{equation}
 0 \leq {< J(x - \bar  x), \bar  x -\xi > }
 =  {< J(x - \bar  x), \bar  x - x > } +  {< J(x - \bar  x),  x -\xi > }
=
\end{equation}
$$ = -||x - \bar  x||^2  +   < {J(x - \bar  x),  x -\xi > }$$
which yields (\ref {f52}).

{\bf  5.f.} Now we describe the property of uniform continuity of
operator $ P_{\Omega}x $
in Banach space $B$. Recall that in  Banach space  the metric projection
operator  is not nonexpansive in general case. But it is uniformly
continuous on
each bounded set according to the following theorem.
\begin {th} \label{f54}
 Let B be an uniformly convex and  uniformly  smooth Banach space.  If
$\delta _B (\epsilon) $  is a modulus of convexity of the space B,
$g_B (\epsilon)  = \delta  _B (\epsilon)  / \epsilon$  and  $g^{-1}_B
(\cdot)$
is an inverse function, then
\begin{equation} \label{f55}
||\bar x - \bar y|| \leq  C  g^{-1}_B  (2LC^2 g^{-1}_{B^{*}} {(2CL||x -
y||)}),
\end{equation}
where $1 < L  < 3.18$ is Figiel's constant (see \cite{aln9}) and
$$ C = 2  {max} \lbrace 1, ||x - \bar y ||, ||y - \bar x||\rbrace .$$
\end {th}
\begin{re}  If  $ ||x - \bar y || \leq R$ and $||y - \bar x|| \leq R$,
then
 $(C = 2 {max} \lbrace 1, R \rbrace) $ is an absolute constant and
(\ref{f55})
provides a quantitative description of the uniform continuity  of
operator $ P_{\Omega} $
in Banach space on each bounded set.
\end{re}
The estimate (\ref{f55}) which was established in \cite{aln2} is global
in nature.
Earlier, in \cite{bj}  Bjernestal obtained local estimate
\begin {equation} \label{f57}
 ||\bar x - \bar y|| \leq 2\delta^{-1}_B (6\rho _B {(2||x - y||)}),
\end{equation}
 where $\rho _B (\tau)$ is  a modulus of smoothness of the space B
\cite{dis} .

The estimate of (\ref{f57}) is better than  our estimate (\ref{f55}).
This is why in
 \cite{al2} we continued the investigation of uniform continuity of
metric
projection operator in Banach space.  It turns out that the following
global variant
of (\ref{f57}) can be obtained.
\begin {th} \label{f58}
 Let B be an uniformly convex and  uniformly  smooth Banach space. If
$\delta _B (\epsilon) $  is a modulus of convexity of the space $B$ and
 $\rho _B (\tau)$ is  a modulus of  its  smoothness, then
\begin {equation} \label{f59}
||\bar x -\bar y|| \leq C\delta^{-1}_B (\rho _B {(8CL||x - y||)}),
\end{equation}
where  constant L and  function C are defined in Theorem \ref{f54}.
\end {th}

\begin{re}  To accuracy constants  the estimates (\ref{f55}) and
(\ref{f59}) give
respectively
 $$||\bar x - \bar y|| \leq  g^{-1}_B  ( g^{-1}_{B^{*}} {(||x - y||)}),
$$
 $$||\bar x - \bar y|| \leq \delta^{-1}_B (\rho _B {(||x - y||)}). $$
\end{re}

{\bf  5.g.} -

{\bf  5.h.} -

{\bf  5.i.} Any  $P_{\Omega}$ satisfies the inequality (see \cite{li})
$$ <J(x - P_{\Omega}x) - J(y - P_{\Omega}y), P_{\Omega}x - P_{\Omega}y> \ge 0,
\;\; \forall x,y \in B.$$

Using  the properties of  metric projection operator $P_\Omega$
we obtained
Banach space analogue of Theorem \ref{th1}.
\begin {th} \label{f280}
 Let $ A $ be an arbitrary operator from Banach space $B $ to $B^{*}$,
$\alpha$ an
 arbitrary fixed positive number, $f\in B^{*}$.  Then  the point $x \in
\Omega  \subset B$
is a solution of variational inequality
 $$<Ax-f, \xi - x> \ge 0, {\qquad }        \forall \xi \in\Omega $$
if and only if x is a solution of the operator equation in B
\begin{equation} \label {f282}
x=P_{\Omega}(x - \alpha J^{*}(Ax-f))
\end{equation}
where  $ J^{*}: B^{\ast} \rightarrow B $ is normalized duality mapping
in $  B^{\ast}$.
\end{th}

Iterative process corresponding to (\ref{f282}) is the following
\begin{equation} \label {f2823}
x_{n + 1} = P_{\Omega}(x_{n}-\alpha_n J^{*} (Ax_{n} - f)),
\qquad n=0,1,2...
\end{equation}
However, there are not any approaches to the investigation of (\ref{f2823}).

In Section 7 and Section 8  we will construct the generalized projection
operato
Banach spaces with the additional properties  5.b,  5.d,  5.g and  5.h.,
and iterative methods for which one can establish convergence, stability,
and nonasymptotic estimates of convergence rate.

\section{Parallelogram Inequalities and Strongly Unique Best
Approximations}
\setcounter{equation}{0}
In this Section we discuss the existence of strongly unique best
approximation in
the spaces $l^p, L^p $  and $ W^p_m$  where  $\infty > p >1. $

In \cite{aln8} (see also \cite{aln9}) we established the following upper
parallelogram inequality
 $$2||x||^2 + 2||y||^2 - ||x + y||^2 \leq 4 ||x - y||^2 +
C_1(||x||,||y||)
\rho _B (||x - y||) ,$$
$$ C_1(||x||,||y||) = 2 \max \lbrace  L, (||x|| +||y||)/2 \rbrace $$
for the arbitrary points $x$ and $y$ from uniformly smooth Banach space
$B$.
We also obtained lower  parallelogram inequality
\begin {equation} \label{f167}
2||x||^2 + 2||y||^2 - ||x + y||^2 \ge L^{-1}\delta_B (||x -
y||/C_2(||x||,||y||)),
\end{equation}
$$ C_2(||x||,||y||) = 2 \max \lbrace  1, \sqrt{(||x||^2
+||y||^2)/2}\rbrace $$
for the arbitrary points   $x$ and $y$ from uniformly convex Banach
space $B$.
Analogous parallelogram inequalities for the $||x||^q$ of other orders
$q$
were obtained in \cite{not}.

If  $||x|| \leq R $ and $||y||\leq R $ then
$$ C_1(||x||,||y||) =
 C_1 = 2  \max \lbrace  L, R\rbrace $$
and
$$C_2(||x||,||y||) = C_2 =
2 \max \lbrace  1, R\rbrace. $$
In this case (\ref {f167}) expresses
the  uniform convexity of  the functional $||x||^2$ on each bounded set
in
$B$ with modulus of convexity
$\delta (||x - y||) = (2L)^{-1}\delta_B (||x - y||/C_2) $, and
\begin {equation} \label{f109}
||x + y||^2  \leq   2||x||^2 + 2||y||^2  - L^{-1}\delta_B (||x -
y||/C_2).
\end{equation}
In \cite{vnc} the following lemma was proved.
\begin {lem} \label{f77}
If a convex functional ${\phi(x)}$ defined on convex closed set $ \Omega $
satisfies the inequality
$${\phi({\frac 12}x + {\frac 12} y)} \leq {\frac 12} \phi (x) +
{\frac 12} \phi (y) -  \kappa (||x - y||),$$
where $\kappa (t) \ge 0,  \kappa (t_0) > 0 $ for some $t_0 > 0$,
then  $\phi (x)$ is uniformly convex functional with modulus of convexity
$\delta (t) = 2 \kappa (t)$ and
$${\phi(y)} \ge  \phi (x) + <l(x),y - x> + 2 \kappa (||x - y||).$$
for all $l(x) \in \partial \phi(x)$. Here $ \partial \phi(x)$ is the
set of all support functionals (the set of all subgradients) of $\phi(x)$ at
the point $x \in \Omega.$
\end{lem}

{}From this lemma and (\ref {f109}) it follows that
\begin {equation} \label{f168}
||x||^2  \leq ||y||^2 +  2<Jx, x-y>  -  (2L)^{-1}\delta_B (||x -
y||/C_2).
\end{equation}
Let  $\Omega \subset B,\; \xi \in \Omega,\; \bar x = P_\Omega x $.  We
replace in (\ref {f168})  $x$ by $(x - \bar x)$ and $y$ by $(x - \xi)$
and
obtain
$$||x - \bar x||^2  \leq ||x - \xi||^2 - 2 <J(x - \bar x), \bar x - \xi>
-  (2L)^{-1}\delta_B (||\bar x - \xi||/C_2).$$
The property  5.c then yields the following general formula
\begin {equation} \label{f169}
 ||x - \bar x||^2  \leq ||x - \xi||^2  -  \lambda \delta_B (||\bar x -
\xi||/C_2), \;\; \lambda = (2L)^{-1}.
\end{equation}
It is obvious  that  if  $\delta_B (\epsilon)$
can be estimated by $\epsilon ^2$
(this occures in Hilbert spaces and in the spaces
of type $L^p$ for  $1 < p \leq 2 $  \cite{aln1}),
then projection $\bar x$ is a strongly
unique best approximation with $\phi (t) = t^2$, at least, on each
bounded set (See Def. \ref {def1}). However, constant  $\lambda $
in inequality (\ref {f169}) is not exact in these cases.

In \cite{Yu} it was shown that in spaces   $W^p _m$ (consequently, in
$L^p $
and  $l^p$),  $1 < p \leq 2 ,$ the following inequality holds
\begin {equation} \label{g69}
 ||x + y||^2  \leq   2||x||^2 + 2||y||^2  - (p-1)||x - y||^2.
\end{equation}
Then Lemma \ref {f77}, the property  5.c and (\ref{g69})  give the estimate
$$||x - \bar x||^2  \leq ||x - \xi||^2  -  (p-1)||\bar x - \xi||^2 /2. $$

Using the inequality (2.3) from \cite{lhh}, we immediately
obtain from Lemma \ref {f77} that in spaces  $L^p $, $1 < p \leq 2 ,$
the estimate
\begin{equation} \label{f161}
||x - \bar x||^2  \leq ||x - \xi||^2  -  (p-1)||\bar x - \xi||^2.
\end{equation}
is valid. This coincides with the result of \cite{lhh} (Theorem 4.1).

Furthermore, the strongly unique best approximation of the projection
$\bar x $  in spaces $B$ of the type  $l^p,\; L^p $  and $ W^p_m$  where
$\infty > p > 2, $
can be established from Lemma \ref {f77} and Clarkson's  inequality
$$ ||x+y||^p  \leq    2^{p-1} ||x||^p  +  2^{p-1}  ||y||^p  - ||x-y||^p.$$
It means that functional $||x||^p$ is uniformly convex.
Therefore, we can write
\begin{equation} \label{f162}
 ||y||^p \ge ||x||^p +p <J^{\mu}x,y-x> + 2^{-p+1} ||x-y||^p
\end{equation}
where  $ J^{\mu}$  is a duality  mapping  with the gauge function $\mu (t) =
t^{p-1}$ (see  Section 7).
Now we substitute in  (\ref{f162}) $(x -\bar x)$ and $(x - \xi)$ for $x$
and $y$,
respectively.  By virtue  of
 $$<  J^{\mu} (x - \bar x), \bar x - \xi > \ge 0 $$
and
$$||x - \bar x||^p \leq  ||x-\xi||^p - p <J^{\mu}(x - \bar x),\bar x - \xi>
- 2^{-p+1}
||\bar x-\xi||^p ,$$
we have
\begin{equation} \label{f177}
||x - \bar x||^p \leq  ||x-\xi||^p  -    2^{-p+1} ||\bar x-\xi||^p .
\end{equation}
This expression improves the estimate of a strongly unique best
approximation.
(Compare with the corresponding inequality in \cite{sm,ls}).
It was obtained without any additional conditions on the modulus of
convexity
of the spaces  and on the sets $\Omega$.
Besides, the important generalization of (\ref {f177}) is valid.
\begin {th}
Let $B$ be a space either $l^p$ or $ L^p $ or $ W^p_m$  where
$\infty > p \ge 2. $ Let $\Omega$ be a closed convex set in $B$. Then for
every point $x \in B$ there exists a unique point $\bar x = P_\Omega x $
such that
$$||x - \bar x||^s \leq  ||x-\xi||^s  -    2^{-s+1} ||\bar x-\xi||^s ,
{\quad} \forall \xi \in\Omega, {\quad} s \ge p \ge 2. $$
\end {th}

The proof follows from the inequality (see \cite {vnc})
 $$||x+y||^s  \leq    2^{s-1} ||x||^s  +  2^{s-1}  ||y||^s  - ||x-y||^s,
  {\quad} s \ge p \ge 2. $$

\section{Generalized Projection Operator $\Pi_\Omega$ in Banach Space}
\setcounter{equation}{0}

Here we introduce generalized projection operator $\Pi_\Omega$  and
describe
its properties in Banach spaces. Then we formulate theorem about
convergence
of the method of succesive projections given in Section 3.  This method
yields
an approximation of the  common  point of convex closed sets. (See
second problem in Section 3).

The  formula (\ref{f1}) in the Definition \ref{gpo} of the metric
projection operators     is equivalent to the minimization problem
\begin{equation} \label{f61}
P_{\Omega}x = \bar x; \;\;\; \bar  x: ||x - \bar x||^2 =
\inf_{\xi \in \Omega} ||x - \xi||^2 .
\end{equation}
Note that $V_1(x, \xi) = ||x - \xi||^2$ can be considered not only as
square
of distance between points $x$ and $\xi$ but also as Lyapunov functional
with respect to $\xi$ with fixed $x$. Therefore,  we can rewrite (\ref
{f61}) in the
form
 $$P_{\Omega}x = \bar x; \;\;\; \bar  x:  V_1(x ,\bar x) =
\inf _{\xi \in \Omega} V_1(x ,\xi). {\qquad } $$
In Hilbert space
 $$V_1(x ,\xi) = ||x||^2 _H-2 (x, \xi) +||\xi||^2_H .$$
In the papers \cite{aln8,aln1} we introduced  Lyapunov functional
\begin{equation} \label{k64}
V_2(Jx ,\xi) = ||Jx||^2_{B^{*}} -2 <Jx, \xi> +||\xi||^2.
\end{equation}
It is a nonstandard functional because it is defined on both the
elements $\xi$ from
 the primary space $B$ and on the elements $(Jx)$ from the dual space
$B^{*}$
(see also (\ref {f72})).

\begin {lem} \label{s1}
The functional $V_2(Jx ,\xi)$ has the following properties:\\
\\
1.  $V_2$ is continuous.\\
\\
2.  $V_2$ is convex with respect to $\varphi = Jx$ when $\xi$ is fixed end
with respect to $\xi$ when $x$ is fixed.\\
\\
3.  $V_2$ is differentiable with respect to $\varphi $  and  $\xi$ .\\
\\
4.  $ grad_{\varphi} V_2(Jx, \xi) = 2(x-\xi).$\\
\\
5.  $ grad_{\xi} V_2(Jx, \xi) = 2(J\xi-Jx).$\\
\\
6.  $V_2(Jx, \xi) > 0, \qquad \forall x,\xi \in B.$\\
\\
7.  $V_2(Jx, \xi) = 0,$ only if  $ x=\xi .$\\
\\
8.  $V_2(Jx, \xi) \rightarrow \infty,$  if  $ ||x||\rightarrow\infty$
(or
$||\xi||\rightarrow\infty$) and vice versa.\\
\\
9.  $(||x||-||\xi||)^2 \leq V_2(Jx, \xi) \leq (||x||+||\xi||)^2 .$\\
\\
10. $L^{-1} \delta_B (||x-\xi|| / C) \leq V_2(Jx, \xi)\leq
L^{-1} \rho_B (8LC||x-\xi||)$, where the constant $L$ is from Theorem
\ref{f54} and
$$C=2 \max \lbrace  1, \sqrt{(||x||^2 +||\xi||^2)/2} \rbrace.$$
11. $V_2(Jx, \xi) \rightarrow 0 $ if $||x - \xi||\rightarrow 0 $, and $||x||,$
 $ ||\xi||$ are bounded, and vice versa .\\
\end{lem}

There  is a connection between the functional $V_2(Jx, \xi)$ and the
Young-Fenchel
transformation, because
 $$||Jx||_{B^{*}} = \sup_{\xi\in B} \lbrace 2<J x,\xi>  -
||\xi||^2 \rbrace .$$
Notice  also that (\ref {k64}) is equivalent to
 $$V_2(Jx ,\xi) = ||x||^2 -2 <Jx, \xi> +||\xi||^2.$$
However, previous form (\ref {k64}) is necessarily used to obtain
properties 4 and 10.
\begin {de}
  Operator $\Pi_\Omega: B \rightarrow \Omega \subset B$ is called
the generalized projection operator if
it puts the arbitrary fixed point $ x \in B$ into the correspondence
with
the point of minimum for the functional $V_2(Jx, \xi)$ according
to the minimization problem
 $$\Pi_{\Omega}x = \hat x; \;\;\; \hat  x:  V_2(Jx ,\hat x) =
\inf _{\xi \in \Omega} V_2(Jx ,\xi) .$$
\end{de}
\begin{re}
 In Hilbert space $ V_2(Jx ,\xi) = V_1(x ,\xi) $ and $\hat x =\bar x$.
 \end{re}
Next we describe the properties of the operator $ \Pi_\Omega$:\\

{\bf  7.a.} The operator $\Pi_{\Omega}$ is fixed in each point $\xi \in
\Omega $ , i.e.
 $\Pi_{\Omega}\xi = \xi .$

{\bf  7.b.}  $\Pi_{\Omega}$ is d-accretive in $B$ , i.e. for all
$x,y \in B $
$$< Jx - Jy, \hat x -\hat y > \ge 0  .$$

{\bf  7.c.}  $ < Jx - J\hat x, \hat x -\xi > \ge 0,  {\qquad }
\forall \xi \in\Omega  .$  \\

{\bf  7.d.}  $ < Jx - J\xi, \hat x -\xi > \ge 0,  {\qquad }
\forall \xi \in\Omega .$ \\

{\bf  7.e. }  $ < Jx - J\hat x ,  x -\xi > \ge 0,  {\qquad }
\forall \xi \in\Omega .$ \\

{\bf  7.f.} $||\hat x -\hat y|| \leq C g^{-1}_B (2LC^2 g^{-1}_{B^{*}}
(2LC ||x-y||))$,
where constant $L$ is  from  Theorem \ref{f54} and
$$ C = 2 \max \lbrace  1, ||x||, ||y||, ||\hat x||, ||\hat y|| \rbrace .$$

\begin{re}   If  $||x ||\leq R, || \hat x || \leq R,  ||y ||\leq R $ and
 $|| \hat y|| \leq R,$ then  $C = 2 {max} \lbrace 1, R \rbrace$  is
absolute constant and 7.f expresses the uniform continuity  of operator
 $\Pi_{\Omega} $ in Banach space on each bounded set.
\end{re}

{\bf  7.g.}   $ < Jx - Jy, \hat x -\hat y > \ge (2L )^{-1}\delta_B
(||\hat x -\hat y||/C),$ where
$$ C = 2 \max \lbrace  1, ||\hat x||, ||\hat y|| \rbrace .$$

{\bf  7.h.}  The operator $ \Pi_{\Omega}$ gives absolutely best
approximation of
$x\in B$ with respect to
 functional $V_2(Jx, \xi)$
 $$V_2(J \hat x, \xi) \leq V_2(Jx, \xi) - V_2(Jx, \hat x). $$
Consequently, $ \Pi_{\Omega}$ is {\it conditionally nonexpanseve \/}
operator in Banach space, i.e.
 $$V_2(J \hat x, \xi) \leq V_2(Jx, \xi).$$

{\bf  7.i.} Any  $\Pi_{\Omega}$ satisfies the inequality
$$ <(J - J\Pi_{\Omega})x - (J - J\Pi_{\Omega})y,
\Pi_{\Omega}x - \Pi_{\Omega}y> \ge 0, \;\; \forall x,y \in B.$$

\begin{re}  If B=H, then the formulas 7.a - 7.e and 7.h - 7.i coincide with
ones
4.a - 4.e and 4.h - 4.i, but 7.f  and 7.g differ from 4.f and 4.g by only
constants
(on any bounded set).
 \end{re}
 Using properties of  the generalized projection operator $\Pi_\Omega$
we
obtained the  theorem.
\begin {th} \label{g4}
The following statements hold for the method of successive generalized
projections
(\ref{f16}) and (\ref{f17}):\\
\\
1)  $V_2(Jx _{n+1},\xi) \leq  V_2(Jx _{n},\xi), {\qquad } \forall \xi
\in\Omega_{*} = \bigcap _{i=1}^{m} \Omega_{i} .$\\
\\
2) There exists  a subsequence $\lbrace x_{n_k} \rbrace $ of the
sequence $ \lbrace x_{n}
 \rbrace$ such that  $x_{n_k} \rightarrow x_{*} $  weakly,
 where $ x_{*}\in  \Omega _{*}.$\\
\\
3) If $J$ is sequential weakly continuous operator then $x_{n}
\rightarrow x_{*} $
 weakly.\\
\\
If $\lbrace x_{n} \rbrace $ is an ordered sequence of the elements $x^j_{i},\;
i = 0,1,...,\; j=m,m-1,...,1,$ such that $x^{m}_i = \Pi_{m}x^{1}_{i-1};\;
x^{j}_i = \Pi_{j}x^{j+1}_i,\; j = m-1,m-2,...,2,1;\; x^{1}_{-1} = x_{0},$
then, in addition to 1) - 3),\\
\\
4)  $\sum _{n=0}^{\infty} V_2(Jx _{n}, x_{n+1}) < \infty .$\\
\\
5)  $ ||x_{n} - x_{n+1}|| \rightarrow 0 ,  \;\;$  for  $n \rightarrow
\infty .$\\
\end {th}
In what follows  we discuss statement 3 from the theorem.
We recall that  $F$ is called a sequentially weakly continuous mapping if
from the relation  $ x_{n} \rightarrow x $ (weakly)  it follows  that
$ Fx_{n} \rightarrow Fx $  (also weakly).

Theorem \ref{g4} is valid for dual mapping $ J^{\mu}$ with the gauge
function
$\mu (t) = t^{p-1}$, defined by the relations
$$||J^{\mu} x||_{B^{*}} = ||x||^{p-1},  {\qquad } <J^{\mu} x, x> = ||x||^{p},
{\qquad } J^{\mu} x = grad ||x||^{p}/p . $$
(Notice that normalized dual mapping  corresponds to $p=2$). We set
$$V_3(J^{\mu} x,\xi) = q^{-1}||J^{\mu} x||^q_{B^{*}}- <J^{\mu} x,\xi>  +
p^{-1}||\xi||^p ,{\qquad } p^{-1} + q^{-1} = 1 .$$
The function $ q^{-1}||J^{\mu} x||^q_{B^{*}} $ is conjugate to the function
$ p^{-1}||x||^p $, i.e.
$$q^{-1}||J^{\mu} x||^q_{B^{*}} = \sup_{\xi\in B} \lbrace <J^{\mu} x,\xi>
- p^{-1}||\xi||^p \rbrace .$$
Therefore  $ V_3(J^{\mu} x,\xi) \ge 0 ,{\quad} \forall  x,\xi\in B$ . If  we
now define
$$ \hat {\hat x}: {\qquad } V_3(J^{\mu} x, \hat {\hat x}) =
\inf_{\xi\in\Omega}V_3(J^{\mu} x,\xi) $$
then it can be shown that\\
$$<J^{\mu} x - J^{\mu} \hat{\hat x}, \hat{\hat x} - \xi> \ge 0, {\qquad}
\forall\xi\in\Omega $$
and
$$V_3(J^{\mu} x, \hat{\hat x})\leq V_3(J^{\mu} x,\xi) - V_3((J^{\mu} \hat{\hat
x},\xi). $$
The last inequalities are used mainly in the proof of Theorem \ref{g4}
which will
be given in forthcoming paper.
The property of the uniform continuity of the operator $J^{\mu} x$ can be
obtained using
the results from \cite{not,zr}.
\begin{co}
 In Banach space $l^p, p>1$ the sequence $ x_{n}$ weakly converges
 to\\  $ x_{*}\in  \Omega _{*}= \bigcap _{i=1}^{m} \Omega_{i}.$\\
 \end{co}
This holds because in  $l^p, \infty >p>1$, dual mapping $ J^{\mu}$ with the
gauge
function $\mu (t) = t^p$  is sequential weakly continuous  \cite{br}.
\begin {re}   $V_3(J^{\mu} x, \xi)$ and $V_2(Jx,\xi)$ coinside for $ p=2$
 (up to constant 2).
 \end{re}
It can be shown in a way
similar to the case of  metric projection operator $P_\Omega$ in Banach
space
\cite{aln2} that generalized projection operator $\Pi_\Omega$ is stable
with
respect to peturbation of the set  $\Omega$.

Let  $\Omega_1$ and  $\Omega_2$ be convex closed sets, $x \in B$ and
 $H (\Omega_1, \Omega_2) \leq \sigma ,$ where
$$ H(\Omega_1, \Omega_2) =
\max \lbrace \sup_{z_1\in \Omega_1} \inf _{z_2\in \Omega_2}
 ||z_1 - z_2||, {\quad}\sup_{z_1\in \Omega_2} \inf _{z_2\in
\Omega_1}||z_1 - z_2||\rbrace $$
is a Hausdorff distance between $\Omega_1$ and $\Omega_2.$
Let also $\Pi_{\Omega_1}x = \hat x_1, \Pi_{\Omega_2}x = \hat x_2$.
\begin {th} \label{g6}
If $B$ is uniformly convex  Banach space, $\delta _B (\epsilon) $  is
modulus of the
convexity of  $B$ and $\delta^{-1}_B (\cdot)$ is an inverse function,
then
$$||\hat x _1 -\hat x_2|| \leq C_1\delta^{-1}_B (4LC_2\sigma ) ,$$
$$ C_1 = 2 \max \lbrace 1,||Jx - J\hat x_1||_{B^{*}},
||Jx - J\hat x_2||_{B^{*}} \rbrace,
$$
$$ C_2 = 2 \max \lbrace ||Jx - J\hat x_1||_{B^{*}},
||Jx - J\hat x_2||_{B^{*}} \rbrace.$$
\end {th}

If  $||x|| \leq R , ||\hat x_1|| \leq R $ and $||\hat x_2 || \leq R ,$
then $ C_1$
and $C_2$ are absolute constants, because operator $J$ is bounded in any
Banach
space.

\section{Generalized Projection Operator $\pi_\Omega$ in Banach Space}
\setcounter{equation}{0}
Here we introduce  generalized projection operator $\pi_\Omega$
in Banach space and describe its properties. Then  we use this operator
to establish equivalence between  solution of the variational inequality
in Banach space and solution of the corresponding operator equation.
In other words, we solve first problem described in Section 2.  Finally,
we obtain a link between operators $\pi_\Omega$ and $\Pi_\Omega$
by means of the normalized duality mappings $J$ and $J^{*}$.

We assume that $ \varphi$ is an arbitrary element of the space $ B^{*}$.
\begin {de}
 The generalized projection $ \tilde \varphi$ of the element  $ \varphi$
on the set $\Omega \subset B$ is given by means of a minimization
problem
 $$\pi_{\Omega}\varphi = \tilde \varphi ;\;\;\;    \tilde \varphi:
V_4(\varphi , \tilde \varphi) = \inf _{\xi \in \Omega} V_4(\varphi ,\xi)
$$
where
\begin{equation} \label{f72}
V_4(\varphi ,\xi) = ||\varphi||^2_{B^{*}} -2 <\varphi, \xi>
+||\xi||^2.
\end{equation}
 \end{de}
\begin {re}   In Hilbert space
$ V_4(\varphi ,\xi) = V_3((J^p \hat{\hat x},\xi) =
 V_2(Jx ,\xi) = V_1(x ,\xi) $ and\\
 $\tilde \varphi = \hat{\hat x }= \hat x =\bar x  $ .
 \end{re}
In what follows we list properties of the generalized projection
operator
$\pi_\Omega$ in Banach Space.\\

{\bf  8.a.}  The operator $\pi_{\Omega}$ is J-fixed in each point $\xi \in
\Omega $, i.e.
  $\pi_{\Omega}J\xi = \xi .$

{\bf  8.b.}  $\pi_{\Omega}$ is monotone in $B^{*}$, i.e. for all
$\varphi_1,\varphi_2 \in B^{*}$
 $$<\varphi_1 - \varphi_2, \tilde \varphi_1 - \tilde \varphi_2 > \ge 0 ,
$$

{\bf  8.c.}  $ <\varphi - J \tilde\varphi, \tilde \varphi - \xi > \ge 0,
 {\qquad }        \forall \xi \in\Omega .$

{\bf  8.d.} $ <\varphi - J \xi, \tilde \varphi - \xi > \ge 0,  {\qquad }
\forall \xi \in\Omega .$

{\bf  8.e.} $ <Jx - J\tilde x, x - \xi> \ge 0$, where $ \tilde x =
\pi_{\Omega}Jx,
           {\qquad }        \forall \xi \in\Omega .$\\

{\bf  8.f.}  $||\tilde \varphi_1 - \tilde \varphi_2||
           \leq C g^{-1}_B (2LC ||\varphi_1 - \varphi_2||_{B^{*}})$,
where the constant $L$ is from Theorem \ref{f54} and
$$ C = 2 \max \lbrace  1, ||\tilde \varphi_1||,|| \tilde \varphi_2||
\rbrace .$$

\begin {re} If
$ ||\tilde \varphi _1|| \leq R   $, $||\tilde \varphi _2|| \leq R$ then
 $C = 2 \max \lbrace 1, R \rbrace $ is an absolute constant and  8.f
expresses
uniform continuity  of  the operator $ \pi_{\Omega} $ in Banach space on
each
 bounded set.
  \end{re}

{\bf  8.g.}  $ <\varphi_1 - \varphi_2, \tilde \varphi_1 - \tilde
\varphi_2 >
           \ge (2L )^{-1}\delta_B (||\tilde \varphi_1 - \tilde
\varphi_2||/C) ,$ where $C$ is the constant from 8.f.

{\bf  8.h.} The operator $ \pi_{\Omega}$ gives absolutely best
approximation of $x\in B$
with respect to functional $V_4(\varphi ,\xi) $
$$V_4(J\tilde \varphi,\xi) \leq  V_4(\varphi ,\xi) -  V_4(\varphi ,
\tilde \varphi). $$
Consequently, $ \pi_{\Omega}$ is {\it conditionally nonexpansive \/}
operator in Banach space, i.e.
$$V_4(J\tilde \varphi,\xi) \leq  V_4(\varphi ,\xi) .$$

{\bf  8.i.} Any  $\pi_{\Omega}$ satisfies the inequality
$$ <(I_{B^{*}} - J\pi_{\Omega})\varphi_1 - (I_{B^{*}} - J\pi_{\Omega})
\varphi_2, \pi_{\Omega}\varphi_1 - \pi_{\Omega}\varphi_2>
\ge 0, \;\; \forall \varphi_1, \varphi_2 \in B^{*} ,$$
where $I_{B^{*}}:B^{\ast} \rightarrow B^{\ast}$ is identical
operator in $B^{\ast}.$

Similarly to operator $\Pi_\Omega $, the generalized projection operator
$\pi_\Omega$  in Banach Space is stable with respect to
perturbation of the set  $\Omega$.

Using  the properties of generalized projection operator $\pi_\Omega$
we obtained
 Banach space analogue of Theorem \ref{th1}.

\begin {th} \label{f80}
 Let $ A $ be an arbitrary  operator from Banach space $B $ to $B^{*}$,
$\alpha$ an
 arbitrary fixed positive number, $f\in B^{*}$.  Then  the point $x \in
\Omega  \subset B$
is a solution of variational inequality
 $$<Ax-f, \xi - x> \ge 0, {\qquad }        \forall \xi \in\Omega $$
if and only if x is a solution of the operator equation in B
\begin{equation} \label {f482}
x=\pi_{\Omega}(Jx - \alpha (Ax-f)).
\end{equation}
\end{th}

It is not hard to verify that
 $$\Pi _\Omega = \pi _\Omega J ,   {\qquad }   \pi _\Omega = \Pi _\Omega
J^{*} $$
 where $J: B\rightarrow B^{\ast}$ is a normalized duality mapping in $B$
 and $ J^{*}: B^{\ast} \rightarrow B $ is normalized duality mapping
in $  B^{\ast}$.  Therefore  we can rewrite (\ref{f482}) in the form of
$$x=\Pi_{\Omega}J^{*}(Jx - \alpha (Ax-f)).$$
Denote also that $ J^{*} = J^{-1}. $

For the iterative process (\ref{k6}) we proved the assertions analogous
to the
Theorems 3 and 4 from the paper \cite{aln1} (see also Remark  7  in
\cite{aln1}).

It is interesting to note that $\pi_\Omega = J^{*}$ in the case $\Omega
= B$
(the problem of solving the equation $Ax = f).$ Then (\ref{f482}) is
rewritten as
$$x=J^{*}(Jx - \alpha (Ax-f)) $$
or in the form of
\begin{equation} \label {f84}
 Jx = Jx - \alpha (Ax-f)
\end{equation}
because $JJ^{*} = I$. Here $I$ is identical operator. Iterative method
for (\ref{f84})
$$Jx_{n+1} = Jx_{n} - \alpha_n (Ax _{n}-f),{\quad }     n=0,1,2,...,
x_o \in B, {\quad }  x_n=J^{*}Jx_n $$
has been studied earlier in  \cite{aln1,aln8}.

Along with (\ref{k6}) we considered the following iterative processes:
$$x_{n + 1} = \pi_{\Omega}(Jx_{n}-\alpha_n (Ax_{n} - f) /||Ax_{n} -
f||),
{\qquad }n=0,1,2... $$
for the variational inequality (\ref{f5}) with nonsmooth unbounded
operator $A,$
and
$$x_{n + 1} = \pi_{\Omega}(Jx_{n}-\alpha_n ( \partial
u(x_{n})/||\partial u(x_{n})||),
{\qquad }n=0,1,2... $$
and
$$x_{n + 1} = \pi_{\Omega}(Jx_{n}-\alpha_n ( u(x_{n}) - u^{*})\partial
u(x_{n})/
||\partial u(x_{n})||^2),
{\qquad }n=0,1,2... $$
for the minimization of functional $u(x).$ Here $u^{*} =
\min_{x  \in  \Omega} u(x).$

\section{Variational Inequalities and Wiener-Hopf Equations in Banach Spaces}
\setcounter{equation}{0}

In  Section 5 and Section 8 we formulated new equivalence theorems between
varia
inequalities in Banach spaces and corresponding operator equations
(\ref{f282}) and (\ref{f482}) with
metric projection operator $P_{\Omega}$  and  generalized projection
operator $\pi_\Omega$, respectively. It is natural to call the equations
of this type $direct$ $projection$ $equations$.
Now we establish the connection of
variational inequalities with other  operator equations, so called
Wiener-Hopf equations.

Let  $P_{\Omega}$ be a metric projection operator in Hilbert space $H,$
$I$ an identical operator, $A$ and $f$ operator and "right hand part" of
variational inequality (\ref {f2}),
$Q_{\Omega} = I-P_{\Omega}$. Then the equation $AP_{\Omega}z + Q_{\Omega}z =f$
is said to be a generalized Wiener-Hopf equations in Hilbert space.

The following theorem is valid in Banach spaces (cf. with Hilbert
case in \cite{shi})
\begin {th} \label{f81}
The variational inequality (\ref {f5}) has a (unique) solution $x \in B $
if and only if the Wiener-Hopf equation
\begin{equation} \label {f99}
J^{*} (AP_{\Omega}z - f) + \alpha^{-1}Q_{\Omega}z = 0
\end{equation}
with an arbitrary fixed positive
parameter $\alpha$ has a (unique) solution $z \in B$ for each $f\in B^{*}.$
Moreover, $z = x - \alpha J^{*}(Ax - f)$ and $x = P_{\Omega}z.$
\end{th}

The simplest iterative method to approximate a solution of the equation
(\ref {f99}) is the following
$$x_n = P_{\Omega}z_n$$
and
$$z_{n+1} =  x_{n} - \alpha_n J^{*} (Ax_{n} - f).$$
However, its convergence can not be established
because of the reasons mentioned above. As before, we will
obtain the convergent iterative process using generalized projection
operators $\pi_\Omega$. But at the beginning we need an equivalence
theorem with this projection operator.
\begin {th} \label{f82}
The variational inequality (\ref {f5}) has a (unique) solution $x \in B$
if and only if the Wiener-Hopf equation
\begin{equation} \label {f100}
A\pi_{\Omega} z + \alpha^{-1} Q_{\Omega} z = f,{\qquad }Q_{\Omega} =
I_{B^{*}} - J \pi_{\Omega},
\end{equation}
with an arbitrary fixed  positive
parameter $\alpha$ has a (unique) solution $z \in B^{*}$ for
each $f \in B^{*}.$
Moreover, $z = J x - \alpha (Ax - f)$ and $x = \pi_{\Omega}z.$
In (\ref {f100}) $I_{B^{*}}:B^{\ast} \rightarrow B^{\ast}$ is identical
operator in $B^{\ast}.$
\end{th}

It turns out that the iterative process
$$x_n = \pi_{\Omega}z_n$$
and
$$z_{n+1} = J x_{n} - \alpha_n (Ax_{n} - f)$$
converges strongly to the solution $z \in B^{*}$ of the Wiener-Hopf equation
(\ref{f100}).  On the other hand, the iterative process
$$x_{n+1} = \pi_{\Omega}z_{n+1}$$
and
$$z_{n+1} = J x_{n} - \alpha_n (Ax_{n} - f)$$
for an approximation of the solution $x$ of the variational inequality
(\ref {f5}) coincides with  (\ref{k6}) and converges strongly too.
\begin {re}
It follows from Theorems \ref{f81} and \ref{f82} that
$$ P_{\Omega}(x - \alpha J^{*}(Ax-f))=\Pi_{\Omega}J^{*}(Jx - \alpha (Ax-f))$$
where $x$ is unique solution of the variational inequality (\ref {f5}).
\end{re}

Similar results can be formulated also for the quasivariational
inequalities and for the complementarity problems. (See, for example,
\cite{jy,no}).

In conclusion , we notice  that generalized projection
operators in Banach spaces constructed above and metric projection
operators in
Hilbert space are defined in similar way in the form of minimization
problems, and these problems are of the same level of difficulty.
This means that generalized projection operators in Banach spaces can be
viewed
as a natural generalization of the metric projection operators in
Hilbert spaces .\\

Results of this paper were  presented at the 1993 SIAM annual meeting
in Philadelphia, July  12-16.

\end{document}